\documentclass{article}
\usepackage{amsfonts}% for mathbb
\usepackage{amsmath}% for align*
\usepackage{graphicx}% for eps documents
\usepackage{amsthm}% for Definition
\newtheorem{theorem}{Theorem}[section]
\newtheorem{proposition}{Proposition}[section]
\newtheorem{lemma}{Lemma}[section]
\newtheorem{corollary}{Corollary}[section]
\theoremstyle{definition}
\newtheorem{definition}{Definition}[section]
\newtheorem{example}{Example}[section]
\begin{document}

\markboth{David A. Krebes}
%{The Group of Units of the $(n,n)$-tangle Monoid is the $n$-strand Braid Group}
{Units of the String Link Monoids}

%%%%%%%%%%%%%%%%%%%%% Publisher's Area please ignore %%%%%%%%%%%%%%
%\catchline{}{}{}{}{}
%%%%%%%%%%%%%%%%%%%%%%%%%%%%%%%%%%%%%%%%%%%%%%%%%%%%%%%%%%%%%%%%%%%

%\title{The Group of Units of the $(n,n)$-tangle Monoid is the $n$-strand Braid Group}
\title{Units of the String Link Monoids}
\author{DAVID A. KREBES \\ dkrebes@gmail.com}

%referee:  Say something about inverses in the whole category of tangles.
%indenting issues
%\address{dkrebes@gmail.com}

\maketitle

\begin{abstract}
We show that the map obtained by viewing a geometric (ie. representative) braid as a string link induces an isomorphism of the $n$-strand braid group onto the group of units of the $n$-strand string link monoid.
\end{abstract}

%\keywords{tangle; braid; unit; string link}

%\ccode{Mathematics Subject Classification 2000: 57M25, 57M27}

\section{Introduction}

In this paper we explore from the algebraic point of view the simple fact that a geometric (ie. member of an equivalence class) braid by definition represents a string link.

The main result of this paper is the following

\begin{theorem}
\label{main}
For each positive integer $n$, the $n$-strand braid group is naturally isomorphic to the group of units of the $n$-strand string link monoid.  Here by ``natural'' we mean induced by the above geometric identification.
\end{theorem}

A proof of this result for pure braids and ``pure'' string links (in the sense of ``pure braids''; see below) can be found in \cite{HM}.  That proof involves mapping class groups and cobordisms; our proof will be more elementary.

In order to prove this result we must establish two things:

\begin{theorem}
\label{inject}
For each $n$, the braid group on $n$ strands injects naturally into the $n$-strand string link monoid.
\end{theorem}

Theorem~\ref{inject} is proved in a more general setting in \cite{Sk}.

\begin{theorem}
\label{unit-braid}
For each $n$ and with definitions as below, an element of the $n$-strand string link monoid is a unit if and only if it has an $n$-strand braid for a representative.
\end{theorem}

From Theorem~\ref{factors} below and properties of braids we can explicitly describe the inverse of a given element, if it has one, as well as settle the issue of one-sided inverses:

\begin{proposition}
If the string link $T$ has a one-sided inverse then that inverse is represented by $\bar{T}$, the reflection of $T$ across the middle plane $x=1/2$ (see below), and is in fact a two-sided inverse.
\end{proposition}

Should we desire to think about the more general notion of a tangle instead of string link, then by the facts that 1) the $n$-strand string link monoid injects naturally into the $(n,n)$-tangle monoid, and that 2) units of the $(n,n)$-tangle monoid must be string links, we have

\begin{proposition}
For each $n$, the group of units of the $(n,n)$-tangle monoid is naturally isomorphic to the group of units of the $n$-strand string link monoid (which in turn is naturally isomorphic to the $n$-strand braid group).
\end{proposition}

We devote a separate section to each of the proofs of the major results Theorems~\ref{inject} (Section~\ref{inject_section}) and \ref{unit-braid} (Section~\ref{unit}).  

Both Theorems~\ref{inject} and \ref{unit-braid} seem to be known already to experts in the field; however the proofs given here, except as indicated, are new to the author.  All of the definitions appearing in this paper are standard, except perhaps for that of string links, which need not be pure.  All the definitions and results of this paper apply for any positive integer $n$.  We work in the piecewise-linear category.  

\section{Preliminaries}

In order to define string links, we begin with the broader concept of {\it tangle.}

A (``geometric'') $(n,n)$-{\it tangle} is a properly embedded subspace of $B=B^3=I\times I\times I$ (the tangle ``ball'') homeomorphic to $n$ copies of $I$ (the ``strands'') whose $2n$ endpoints are, in no particular order:

$$(0,\frac{1}{n+1},\frac{1}{2}), (0,\frac{2}{n+1},\frac{1}{2}), ...  (0,\frac{n}{n+1},\frac{1}{2}), (1,\frac{1}{n+1},\frac{1}{2}), (1,\frac{2}{n+1},\frac{1}{2}), ...  (1,\frac{n}{n+1},\frac{1}{2})$$
\vskip 1em
%$$ $$
\noindent (Here $I$ is the interval $[0,1]$ of the real line.  In the literature, eg. \cite{Kr}, tangles are sometimes allowed to have loop components but here we disallow them for simplicity.)  This choice of endpoints allows us to compose tangles; see the next section.

In the figures the first coordinate corresponds to lateral position (axis pointing right);  the second to height (axis pointing up); and the third to depth (axis pointing out of page).

%''manifold'' or ``subspace'' or ``subset''?

\begin{definition} For each $n$, the subclass of {\it $n$-strand string links} consists of those $(n,n)$-tangles in which each strand has one endpoint on $\{0\}\times I\times I$ (in the figures, the left face) and the other on $\{1\}\times I\times I$ (the right face).  Note: Contrary for example to \cite{HL} we do not require string links to be pure, in the sense of pure braids (ie. a pure string link is one in which the endpoints of each strand coincide in the second and third coordinates.)
\end{definition}

It is easy to tell whether or not a given tangle is a string link.  For this reason we will restrict ourselves to the class of string links.

\begin{definition}
\label{braid} A {\it braid} is a string link whose strands are realized as images of embeddings from $I$ to the tangle ball $B^3$ in which the first coordinate is a strictly monotone function on $I$.  It follows that the intersection of the braid with a level set $x=\hbox{const}$ consists of $n$ points in the interior of the disk $\{\hbox{const}\}\times I \times I$.   A braid may also be described as a continuous map $f:I\times \Sigma\rightarrow \mathring{D}^2$ whose cross-sections $f_t$, $t\in I$ are injective and such that $f(\{0\}\times \Sigma)$ and $f(\{1\}\times \Sigma)$ occupy the standard positions $(\frac{1}{n+1},\frac{1}{2}), (\frac{2}{n+1},\frac{1}{2}), ...  (\frac{n}{n+1},\frac{1}{2})$ in $D^2$, where $\mathring{D}^2$ is the interior of $D^2=I\times I$ and $\Sigma=\{1,...,n\}$.  The corresponding string link of the first formulation is the union over $t\in I$ and $i\in \Sigma$ of the points $(t, f(t,i))$.

\end{definition}

\section{The Monoid Construction}
\begin{definition} 
\label{tangle_equiv}
Two tangles $T$ and $T'$ are {\it equivalent} if there is a self-homeo-morphism of the tangle ball $B$ that is the identity on the boundary and throws $T$ onto $T'$.  By a theorem of Alexander and Tietze (see \cite{BZ} p. 5) the existence of such a self-homeomorphism implies the existence of an {\it ambient isotopy} $H:B^3\times I\rightarrow B^3$ with the usual properties a) $H_0$ is the identity on $B^3$; b) $H_1$ throws $T$ onto $T'$; and c) $H_t$ is a self-homeomorphism of $B^3$ which point-wise fixes the boundary of $B$ for all $t\in I$.
\end{definition}

As is usual in the literature we shall abuse notation by using the terms ``tangle'', ``string link'' and ``braid'' to refer either to a specific embedded manifold (ie. ``geometric'' object) satisfying the appropriate conditions (which we may speak of as being ``equivalent to'' another such manifold) or an equivalence class of such manifolds (which may be ``represented by'' a particular manifold.)  It will always be clear from the context (eg. by these terms ``equivalent to'' or ``represented by'') which meaning is meant.

\begin{figure}
\begin{center}
\includegraphics[width=5truein]{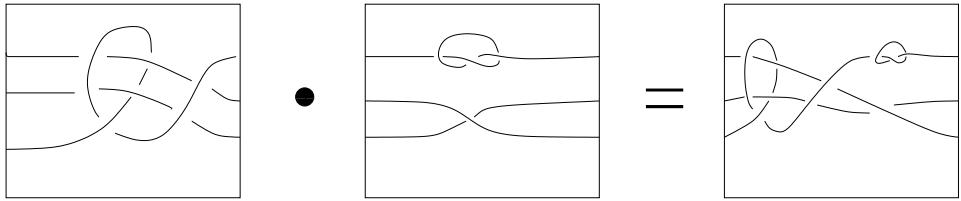}
\caption{An Example of String Link Composition ($n=3$)}
\label{composition}
\end{center}
\end{figure}

\begin{figure}
\begin{center}
\includegraphics[width=2.0truein]{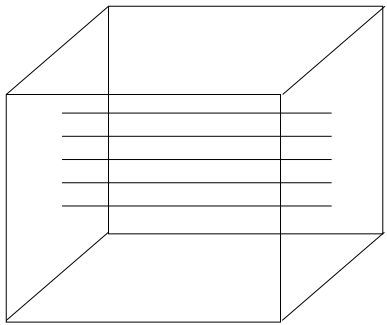}
\caption{The Identity Braid $\beta_e$.  ($n=5$)}
\label{identity}
\end{center}
\end{figure}

The equivalence classes of $(n,n)$-tangles form a monoid (ie. a set with associative binary operation and an identity element, as below) where the binary operation (``composition'' or ``product'') is performed on the representative level by concatenation (and horizontal compression).  See Figure~\ref{composition}, where the tangles happen to be string links.  This operation does not depend on the representatives chosen.  The identity element is represented by $n$ horizontal strands and is denoted $\beta_e$:  See Figure~\ref{identity}.  ($\beta_e$ also happens to be a string link, in fact a braid.)  This same operation (and identity element $\beta_e$) is used in the ``string link monoid'' and the ``braid monoid'' (in fact ``braid group'') for given $n$.

For braids the equivalence relation is defined differently.  We say that $\beta$ and $\beta '$ are {\it isotopic} if there is a continuous family $\beta_t$, $t\in I$ of braids with $\beta_0=\beta$ and $\beta_1=\beta '$.  The elements of the braid group are isotopy classes of braids.  But see Theorem~\ref{equiv} below.
%The operation is well-defined in either case.

Recall that a {\it unit} of a monoid $M$ with identity element $e$ is an element $g\in M$ such that $gh=hg=e$ for some $h\in M$, denoted $g^{-1}$.  The units form a group, the {\it group of units}.  It can be shown for example that neither of the factors of Figure~\ref{composition} is a unit (for the first factor see Section~\ref{examples} below).

%\begin{definition}
%The construction of this section can be applied to the more specific class of $n$-strand string links to yield the $n$-{\it strand string link monoid}, with which we will deal excusively for the rest of this paper.  In particular, the composition of two string links is a string link.
%\end{definition}

\section{Injectivity}
\label{inject_section}
This section is devoted to a proof of Theorem~\ref{inject}, which can be re-formulated as

\begin{theorem}
\label{equiv}
Two braids $\beta$ and $\beta '$ are isotopic as braids if and only if they are equivalent as string links.
\end{theorem}

Theorem~\ref{equiv} implies Theorem~\ref{inject}:  The forward direction shows that what we are calling the ``natural'' map is well-defined and the reverse direction shows that it is injective.

Theorem~\ref{equiv}, forward direction, which is known as the ``Braid Isotopy Extension Theorem'' is proven in \cite{Ka} (Theorem 1.11) and \cite{A} (Theorem 6).

%  This is a special case of an appropriate version of the Isotopy Extension Theorem.

%follows from an appropriate version of the Isotopy Extension Theorem.  Ad hoc proofs appear in \cite{A}, Theorem 6 and \cite{Ka}, Theorem 1.11.

%(and the second formulation of Definition~\ref{braid} above). Note that either proof hands us not just an equivalence but an ambient isotopy, and furthermore one which is level-preserving.
%For an ad hoc proof, the reader is referred to \cite{A}, Theorem 6 (and the second formulation of Definition~\ref{braid} above). Note that either proof hands us not just an equivalence but an ambient isotopy, and furthermore one which is level-preserving.

For the reverse direction, suppose that the braids $\beta$ and $\beta '$ are equivalent as string links.  The existence of the self-homeomorphism of Definition~\ref{tangle_equiv} implies that $\beta$ and $\beta '$ yield the same automorphism $F_n\rightarrow F_n$, where $F_n$ is the fundamental group of the left (right) face of the tangle ball with tangle endpoints removed, a free group of rank $n$.  See \cite{R1} p. 10.  By the faithfulness of this representation (proven by Artin in \cite{A}, Theorem 14) we conclude that $\beta = \beta '$ as braids. 
\hfill \qed
%\hfill$\Box$

\begin{corollary}
\label{inverse}
Two braids $\beta$ and $\beta '$ are inverses as braids if and only if they are inverses as string links.
\end{corollary}

{\it Proof.} Apply Theorem~\ref{equiv} to $\beta\cdot\beta '$ (as well as $\beta '\cdot\beta$) and $\beta_e$. \hfill \qed %\hfill$\Box$ 

\section{A Unit String Link is Represented by a Braid}
\label{unit}
We now turn to the proof of Theorem~\ref{unit-braid}.

Since a braid $\beta$ represents a unit of the string link monoid (by Corollary~\ref{inverse}), we have the reverse direction immediately.

For the forward direction, we will use

\begin{theorem}
\label{factors}
Let $T$ and $T'$ be string links. Then $T\cdot T'$ is equivalent to a braid if and only if both $T$ and $T'$ are equivalent to braids.
\end{theorem}

Setting $T\cdot T'$ to be the braid $\beta_e$ we see that the forward direction of Theorem~\ref{unit-braid} follows.

{\it Proof of Theorem~\ref{factors}.} The reverse direction is immediate from definition.  For the forward direction, suppose that $T\cdot T'$ is equivalent to the braid $\gamma$.  (We use Greek letters for braids.)  We will show below the proof of Lemma~\ref{disks} that $T$ is equivalent to a braid $\beta$.   Then $T'$ is equivalent to the braid $\beta^{-1} \cdot  \gamma$:

\begin{align*}
T\cdot T'&\cong\gamma\\
\beta \cdot T'&\cong\gamma\\
\beta ^{-1}\cdot \beta\cdot T'&\cong\beta^{-1}\cdot \gamma\\
\beta_e\cdot T'&\cong\beta ^{-1}\cdot \gamma\tag{by Corollary~\ref{inverse}}\\
T'&\cong\beta ^{-1}\cdot \gamma\\
\end{align*}

%''represented by'' - element of the monoid.  algebraic.  ``equivalent to'' manifold.  geometric

%Already tried putting ``\hfill'' outside of $$'s.  Okay, but got an extra unwanted line break.
(Here ``$\cong$'' denotes equivalence as string links.)
%T\cdot T'=T_e  b\cdot T'=T_e  b'\cdot b\cdot T'=b'

%We seem to be applying Thm. equiv to b'\cdot b and T_e in both directions.

%\begin{lemma}
%\label{disks}
%referee:  Re-do this using ribbons (Centralisers paper). ``So are these arcs'' fixed 13/3/2013
%A string link $T$ with strands $\alpha_1,...\alpha_n$ is equivalent to a braid if and only if the strands are simultaneously boundary parallel in the following sense:  There are $n$ pairwise disjoint disks $D_1, ...,D_n$ in the tangle ball of which $D_i$ is cobounded by two arcs: $\alpha_i$ and $\beta_i$, where $\beta_i$ is an arc on the boundary of the tangle ball, as shown in Figure~\ref{beta_i_arcs}.  The $\beta_i$'s are, except for the subarcs on the right face of the tangle ball $B^3$, the intersections of the boundary of $B^3$ with $n$ horizontal half-planes $\{(x,\frac{i}{n+1},z)|x,z\in\mathbb{R}; z\ge\frac{1}{2}\}$.  The subarcs on the right face are, except for their endpoints, arbitrary.  The right endpoint of $\beta_i$ of course coincides with the right endpoint of $\alpha_i$.  Since the $D_i$'s are disjoint, so are these arcs $\alpha_i\cup\beta_i$ for different $i$'s.
%\end{lemma}

\begin{figure}
\begin{center}
\includegraphics[width=2in]{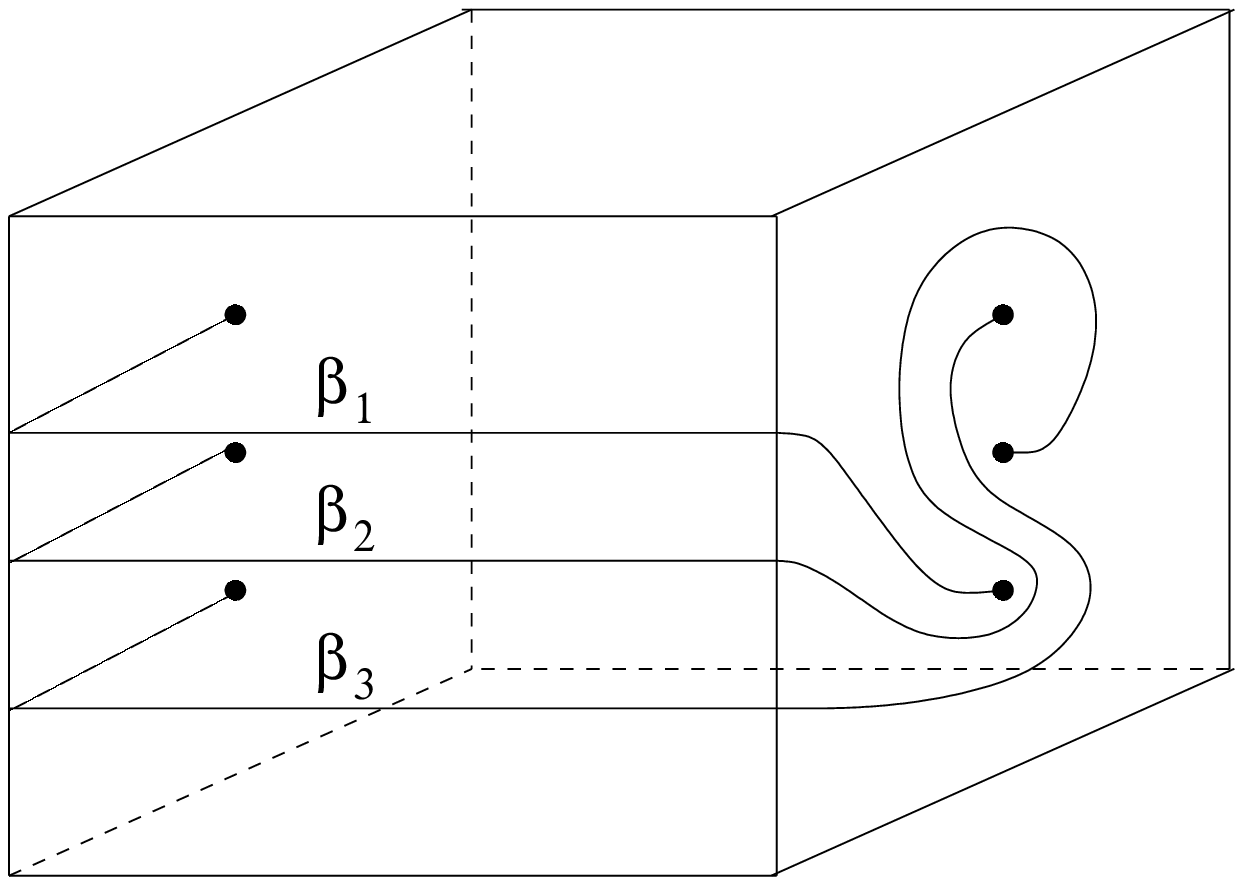}
\caption{The Co-bounding Arcs Construction for Braids ($n=3$).  The subarcs of the $\beta_i$'s on the right face are, except for their endpoints, arbitrary.}
\label{beta_i_arcs}
\end{center}
\end{figure}

We will use a characterization of braids based on the following definition:

\begin{definition}
\label{boundary_parallel}  The strands $\alpha_1,...\alpha_n$, ordered by their left endpoints, of a string link are said to be simultaneously boundary parallel in the sense of this definition if the following situation obtains:  There are $n$ pairwise disjoint disks $D_1, ...,D_n$ in the tangle ball of which $D_i$ is cobounded by two arcs: $\alpha_i$ and $\beta_i$, where $\beta_i$ is an arc on the boundary of the tangle ball. The $\beta_i$'s are, except for the subarcs on the right face of the tangle ball $B^3$, the intersections of the boundary of $B^3$ with $n$ horizontal half-planes $\{(x,\frac{i}{n+1},z)|x,z\in\mathbb{R}; z\ge\frac{1}{2}\}$.  The subarcs on the right face are, except for their endpoints, arbitrary.  The right endpoint of $\beta_i$ of course coincides with the right endpoint of $\alpha_i$.  See Figure~\ref{beta_i_arcs}. 
%Since the $D_i$'s are disjoint, so are these arcs $\beta_i$ for different $i$'s.
\end{definition}

\begin{lemma}
\label{disks}
%referee:  Re-do this using ribbons (Centralisers paper). ``So are these arcs'' fixed 13/3/2013
A string link $T$ is equivalent to a braid if and only if the strands are simultaneously boundary parallel in the sense of Definition~\ref{boundary_parallel}.
\end{lemma}

{\it Proof of Lemma.}

For the forward direction, note that the existence of the disks depends only on the equivalence class of the string link $T$.  Thus we may assume that $T$ is a braid.  An adaptation of the proof of Artin's theorem mentioned above shows that the isotopy $I\times \Sigma\rightarrow D^2$ (see Definition~\ref{braid} of the present paper, second formulation) extends to an isotopy $H:I\times D^2\rightarrow D^2$ point-wise fixing the boundary of $D^2$.  Choose $D_i$ to be the union over $t\in I$ of the arcs $(t,H(t,\frac{i}{n+1},\{z|\frac{1}{2}\le z\le 1\}))\subset B^3$.  (Here $(\frac{i}{n+1},z)\in D^2=I\times I$).
%                        $(t,H(t,\frac{i}{n+1},\{z|\frac{1}{2}\le z\le 1\}))\subset
% changed the following first sentence 13/3/2013
For the reverse direction, we can slide $\alpha_i$ across $D_i$ until it is very close to $\beta_i$ and the (projection onto the) first coordinate pre-composed with any parametrization of the arc is a monotone function, as required for braids. We can easily extend this into a homeomorphism of the surrounding space. Thus $T$ is equivalent to a braid. \hfill \qed %\hfill$\Box$

\vskip 1em
%{\it Proof of Theorem~\ref{unit-braid}, forward direction, continued}
We now return to the proof of Theorem~\ref{factors}.  Specifically, we need to show that the string link $T$ is equivalent to a braid.

Now our assumption that $T\cdot T'$ is equivalent to a braid together with the lemma hands us a set of disks $D_1,...D_n$ with the stated properties.  Let $P$ be the intersection (a disk) of the two tangle balls for $T$ and for $T'$ respectively. We can assume that these $D_i$'s are all transverse to $P$.  Each $D_i$ intersects $P$ in at least a single arc with one endpoint in $\beta_i$ and the other in $\alpha_i$.  (No arc can join points of distinct $D_i$'s because these disks are disjoint.)  This uses up all $n$ points of $(T\cdot T')\cap P$, so that the only other components of intersection of $P$ and $D_1\cup ... \cup D_n$ are loops (simple closed curves).

We wish to arrange that this set of loops is empty.  For this we use the method of ``innermost disks'' as in the proof of the additivity of knot genus (see for example 5A14 in \cite{R}). To do this consider the aggregate (over all of $i=1,...n$) of these loops in $P$ and pick an inner-most one $L$.  This loop $L$ also bounds a disk (which may contain other loops) in $D_i$ for some $i$; replacing it with a parallel copy of the disk bounded by $L$ in $P$ reduces (possibly by more than one) the number of loops of intersection.  Continuing in this way we can remove all of them, as desired.  

Discarding everything to the right of $P$ we still have $T$ sitting inside its tangle ball (with $P$ as the right face) but we also have a new set of disks (remnants of the revised $D_i$'s) inside this ball each of whose boundary is the union of two arcs as in Definition~\ref{boundary_parallel}.  Thus $T$ is equivalent to a braid, as claimed. \hfill \qed %\hfill$\Box$

\section{Special Cases and a Question}
\label{examples}

\begin{example}
\label{integral}
Case $n=2$.  The units of the $2$-strand string link monoid are represented by the braids $\{\sigma_1^k | k\in \mathbb{Z}\}$.  The group of these units is isomorphic to the additive group of integers $\mathbb{Z}$ and we therefore refer to them as the {\it integral} string links (tangles). (See \cite{Kr} p. 342.)
\end{example}

%Proposition.   If we dissolve some strands of a unit we still get a unit.
%\section{Criteria}
%It appears to the author that it can be a difficult task to determine whether or not a given string link is equivalent to a braid.  However we do have the following two results:
%\begin{proposition}
%The string link $T$ is equivalent to a braid if and only if $T \cdot \bar{T} \cong \beta_e$, where $\bar{T}$ is the reflection of $T$ across the plane $x=1/2$.
%\end{proposition}
%The easy proof is omitted.
% by remarking that the author is not aware of any general test of whether or not a given string link is equivalent to a braid.  However, we have the following

We now give an example of something that is not a braid:

\begin{example} Consider a product $T$ of three $n$-strand string links, $n\ge 2$, the second of these in its tangle ball $B$ as depicted in Figure~\ref{not-braid} (shown for large $n$ and vertically expanded for clarity.)  We will show that $T$ is not equivalent to a braid.

\begin{figure}
\begin{center}
\includegraphics[width=3.5in]{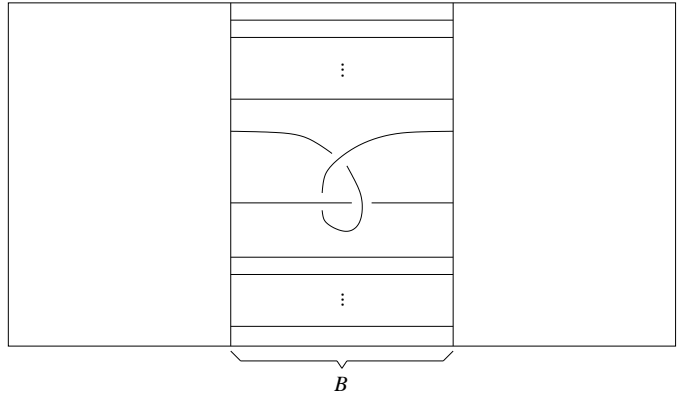}
\caption{A composition of three $n$-strand string links with this middle factor is not equivalent to a braid.}
\label{not-braid}
\end{center}
\end{figure}

Assume then that $T$ is equivalent to a braid.  Then the $n$ arcs of the figure come from distinct strands of $T$.  This is because every strand of $T$ intersects the middle ball $B$ (shown) at least once, and by the pigeon-hole principle, only once.  After dissolving $n-2$ horizontal strands of $T$ leaving the two strands in the middle, and adjusting the endpoints, we are left with a $2$-strand string link $T'\cdot S\cdot T''$ intersecting $B$ in the $2$-strand string link $S$ on the left side of Figure~\ref{trefoil}.  If $T$ is equivalent to a braid then $T'\cdot S\cdot T''$ is equivalent to a braid and by Theorem~\ref{factors} so is $T'\cdot S$ and by another application so is $S$.  This means that $S$ is an integral tangle (see Example~\ref{integral} above) and the closure $d(S)$, shown on the right side of Figure~\ref{trefoil}, is the unknot rather than the trefoil shown.  This contradiction proves the result.\hfill \qed %\hfill$\Box$

\begin{figure}
\begin{center}
\includegraphics[width=2in]{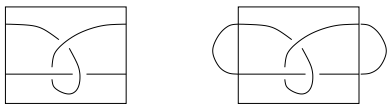}
\caption{The string link $S$, left, and $d(S)$}
\label{trefoil}
\end{center}
\end{figure}

\end{example}

We conclude with the advice of a referee, which we leave for the reader to address:  Say something about (one-sided) inverses in the whole category whose objects are non-negative integers and the morphisms from $k$ to $l$ are the $(k,l)$-tangles.

\section{Acknowledgements}

We wish to thank Thomas Fiedler for pointing the author in the direction of the braid action on the free group used in the proof of Theorem~\ref{equiv}, reverse direction and Mark Grant for the proof of the forward direction of the same theorem in \cite{Ka}.

\end{document}